\documentclass[a4paper,11pt]{article}
\usepackage{longtable}
\usepackage{theorem}
\usepackage[english]{babel}
\usepackage{amsmath,amssymb,amsfonts}
\usepackage{graphics}

\usepackage{graphicx}
\usepackage{latexsym}
\usepackage{float,html}

\newtheorem{moriyamathm}{Theorem}
\newtheorem{prp}[moriyamathm]{Proposition}
\newtheorem{moriyamalem}[moriyamathm]{Lemma}
\newtheorem{cor}[moriyamathm]{Corollary}

\newtheorem{rem}[moriyamathm]{Remark}
\theorembodyfont{\rmfamily}
\newtheorem{moriyamadf}{Definition}

\newcommand{\vect}[1]{{#1}}  
\renewcommand{\AA}{\mathcal{A}}

\newcommand{\DD}{\mathcal{D}}

\newcommand{\LL}{\mathcal{L}}

\newcommand{\R}{{\mathbb{R}}}

\newenvironment{Proof}
{\begin{rm}\par\smallskip\noindent{\bf Proof.}\quad}{\hfill$\Box$\end{rm}}

\begin{document}

\title{%
The Holt-Klee condition for oriented matroids}
\author{
Komei Fukuda%
\thanks{Research partially supported by the Swiss National Science Foundation Project  200021-105202.}\\
Swiss Federal Institute of Technology \\
Zurich and Lausanne, Switzerland\\
fukuda@ifor.math.ethz.ch
\and
Sonoko Moriyama%
\thanks{Research partially supported by Grant-in-Aid for Scientific
    Research from Ministry of Education, Science and Culture, Japan, and
    Japan Society for the Promotion of Science.}\\
Graduate School of Information Science \\
and Technology, \\
University of Tokyo, Japan\\
moriso@is.s.u-tokyo.ac.jp
\and
Yoshio Okamoto%
\thanks{Research partially supported by Grant-in-Aid for Scientific
    Research from Ministry of Education, Science and Culture, Japan, and
    Japan Society for the Promotion of Science.}\\
Department of Information and Computer Sciences,\\
Toyohashi University of Technology, Japan\\
okamotoy@ics.tut.ac.jp}

\date{May 30, 2007}

\maketitle

\begin{abstract}
Holt and Klee have recently shown that every (generic) 
LP orientation of the graph of a $d$-polytope satisfies 
a directed version of the $d$-connectivity property,
i.e.\ there are $d$ internally disjoint directed paths from 
a unique source to a unique sink.
We introduce two new classes HK and HK* of oriented matroids (OMs) by 
enforcing this property and its dual interpretation in terms of line shellings,
respectively.  Both classes contain all representable OMs by the Holt-Klee theorem.
While we give a construction of an infinite family of non-HK* OMs,
it is not clear whether there exists any non-HK OM. 
This leads to a fundamental question as to whether 
the Holt-Klee theorem can be proven combinatorially by using the OM axioms only.
Finally, we give the complete classification of OM(4, 8), the OMs of
rank $4$ on $8$-element ground set with respect to the HK, HK*, Euclidean and Shannon properties.
Our classification shows that there exists no non-HK OM in this class.
\end{abstract}

\section{Introduction}

Let $P$ be a $d$-dimensional convex polytope ($d$-polytope)
in $\R^d$.
We consider a linear program whose feasible region is $P$
with a generic objective function
$f(\vect{x}) = \vect{c}^T \vect{x}$,
i.e. $f(\vect{u}) \neq f(\vect{v})$
for any two distinct vertices $u$ and $v$.
We orient each edge $(u, v)$ from $u$ to $v$ 
if and only if $f(\vect{u}) < f(\vect{v})$.
The resulting orientation on the graph $G(P)$ of $P$
is known as an \emph{LP orientation},
which represents the possible pivot operations
of the simplex method to solve the linear program.
We call $G(P)$ with an LP orientation an \emph{LP digraph}.

\begin{figure}[ht]
\begin{center}
\includegraphics[scale=0.35]{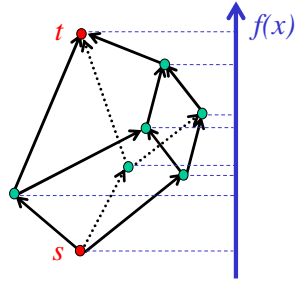}
\end{center}
\caption{An LP digraph induced by $f$}
\label{fig:LPgraph}
\end{figure}

Every LP digraph satisfies the following three properties:
(1) acyclicity, i.e., there exists no directed cycle,
(2) \emph{unique sink\&source property}~\cite{SW01},
namely there exist a unique sink and a unique source,
and (3) the \emph{Holt-Klee condition}~\cite{HK}.
The Holt-Klee condition is a directed version
of the $d$-connectivity property by Balinski \cite{Bal61},
i.e., there are
$d$ internally disjoint directed paths from the source to the sink.
Particularly when $d=3$, these three properties
are also sufficient for an LP digraph~\cite{MK00}.

\begin{figure}[ht]
\begin{center}
\includegraphics[scale=0.40]{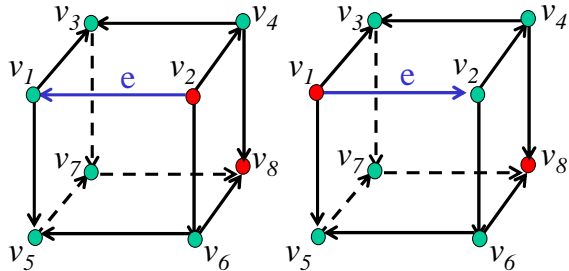}
\end{center}
\caption{Two orientations satisfying and not satisfying the Holt-Klee condition.}
\label{fig:HK-non-HK}
\end{figure}

For example, consider two orientations
on the graph $G(C_3)$ of a $3$-cube $C_3$
in Figure \ref{fig:HK-non-HK}.
While the minimum size of vertex cut sets is three
in the left orientation,
it is only two in the right orientation and in particular,
every dipath from the source $v_1$ to the sink $v_8$ must 
go though either $v_2$ or $v_7$.
Hence, the right digraph does not satisfy the Holt-Klee condition and thus is not an LP digraph.

Our main objective is
to understand how restrictive the Holt-Klee condition is.
The original proof of the Holt-Klee theorem \cite{HK} 
relies heavily on geometric operations
such as affine transformations and orthogonal projections,
and it is not clear whether more combinatorial proof is possible or not.
In particular, it is natural to ask whether this condition is valid
for the oriented matroid program,
which is a combinatorial abstraction of the linear program.
This motivates us to study
the Holt-Klee condition in the setting of \emph{oriented matroids}.
For this purpose,
we introduce two new subclasses of oriented matroids using this condition,
the class of \emph{HK} oriented matroids
and the class of \emph{HK*} oriented matroids.

Before defining the two subclasses,
we set notations for oriented matroids.
We assume that the reader is familiar with oriented matroids.
The standard reference is \cite{om}.
An oriented matroid $M$ is defined as a pair $(E, \LL)$ of
a finite ground set $E$ and the set $\LL$ of covectors.
By the topological representation theorem,
every oriented matroid of rank $d+1$ can be
represented by a pseudosphere arrangement~\cite{FL}.
A \emph{pseudosphere arrangement\/}
 is specified as a triple $\AA=(E,S,\DD)$ of
a finite ground set $E$,
a $d$-dimensional unit sphere $S$ in $\R^{d+1}$, and
a family $\DD = \{ \{ s_e^+, s_e, s_e^- \} : e \in E \}$, 
where $s_e$ is a $(d{-}1)$-dimensional pseudosphere on $S$,
$s_e^+$ is the \emph{positive side\/} of $s_e$, and
$s_e^-$ is the \emph{negative side\/} of $s_e$, respectively.
The \emph{location vector\/} of a point $x \in S$
is the sign vector $\sigma(x) \in \{ +, 0, - \}^E$
defined by $\sigma(x)_e = +$ if $x \in s_e^+$,
$\sigma(x)_e = 0$ if $x \in s_e$,
and $\sigma(x)_e = -$ if $x \in s_e^-$.
A pseudosphere arrangement $\AA$ is said to \emph{represent\/} an oriented matroid
$M = (E, \LL)$ if $\sigma(S) \cup \{ \mathbf{0} \} = \LL$,
where  $\sigma(S) :=\{ \sigma(x) : x \in S \}$.
For the sequel, we use $s_e$
as a topological representation of an element $e \in E$ of $M$.
If all pseudospheres $\{ s_e :e \in E \}$ are realized
by $(d{-}1)$-dimensional linear spheres on $S$,
i.e. the intersection of $S$ and a hyperplane through the origin,
$M$ is said to be \emph{representable},
and \emph{non-representable} otherwise.
We will define the HK property and the HK* property of oriented matroids
to be a consistent generalization of the Holt-Klee condition,
extended to oriented matroids.

First, we define the class of HK oriented matroids,
which is based on the direct relation
between the oriented matroid program and the Holt-Klee condition.
An \emph{oriented matroid program} is a triple $\pi = (M, g, f)$
where $M$ is an oriented matroid $M=(E, \LL)$,
$g\in E$ is not a loop of $M$,
and $f(\neq g) \in E$ is not a coloop of $M$.
The {\em feasible region} $P_\pi$ of $(M, g, f)$ is
$P_\pi := \{ X \in \LL : X_g = +,
X_e \in \{+,0\} \text{ for all } e \in E \setminus \{ g, f \} \}$.
The region $P_\pi$ is \emph{unbounded}
if it is nonempty and there exists $X \in \LL$
such that $X_g = 0$ and $X_e = \{+,0\}$ for all $e \in E \setminus \{ g, f \}$,
and \emph{bounded} otherwise.
We orient the graph of $(M, g ,f)$,
i.e., the $1$-skeleton of the arrangement $\{ s_e : e \in E \setminus \{ f \} \}$
restricted to the positive side $s_g^+$ of the infinity $g$,
so that each edge is oriented from the negative side $s_f^-$ of the objective $f$
toward the positive side $s_f^+$ of $f$
in~\cite[Definition 10.1.16]{om}.
The graph $G_\pi$ restricted to $P_\pi$ is denoted by $G_\pi^+$
and called an \emph{OMP digraph}.
The objective $f$ is \emph{generic} in $P_\pi$
if there is no non-oriented edges in $G_\pi^+$.
An oriented matroid program $\pi = (M, g, f)$
is called \emph{proper} if $P_\pi$ is bounded,
full dimensional (i.e. containing a tope)
 and $f$ is generic.
\begin{figure}[ht]
\begin{center}
\includegraphics[scale=0.38]{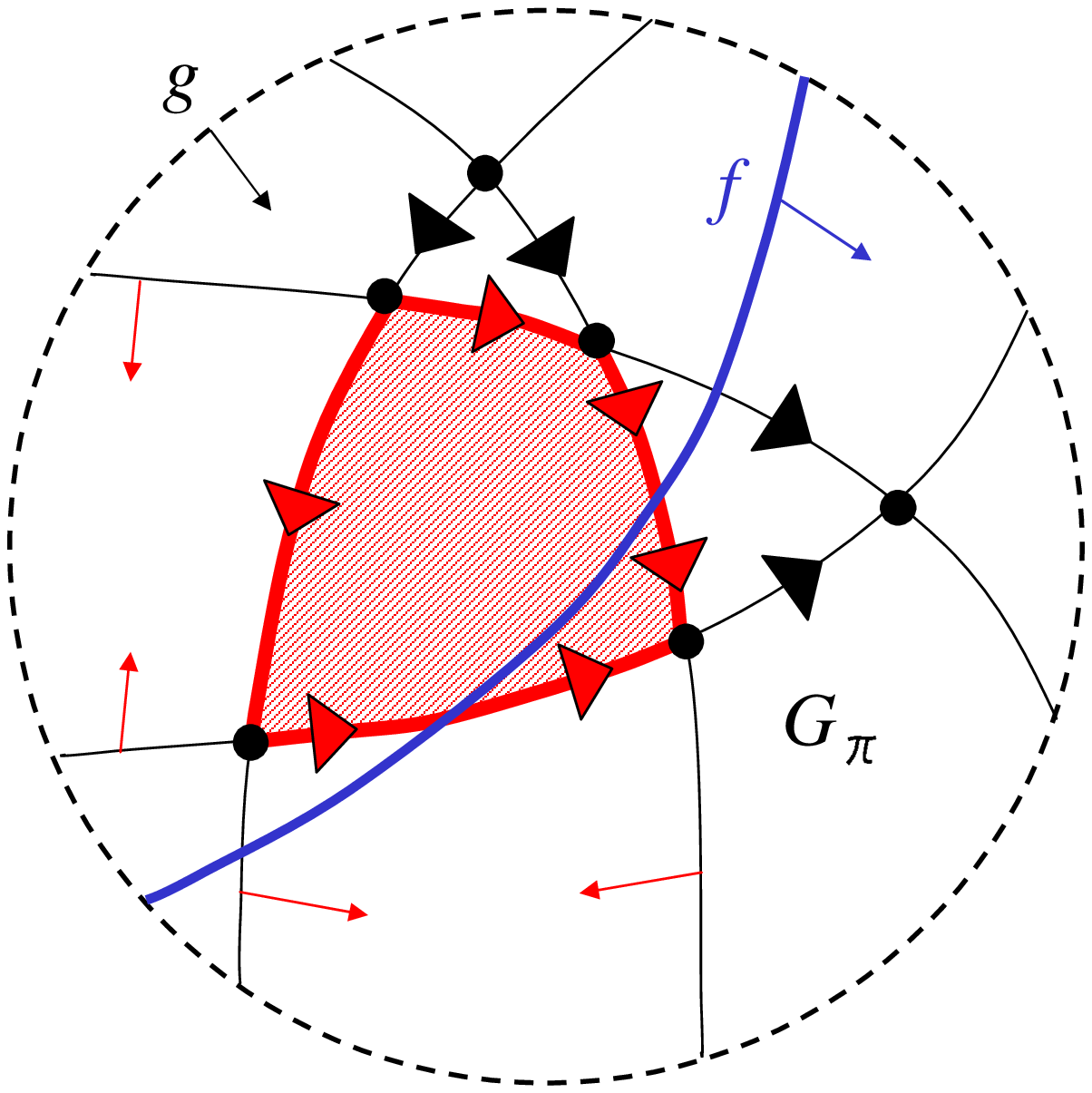}
\hspace{8mm}
\includegraphics[scale=0.38]{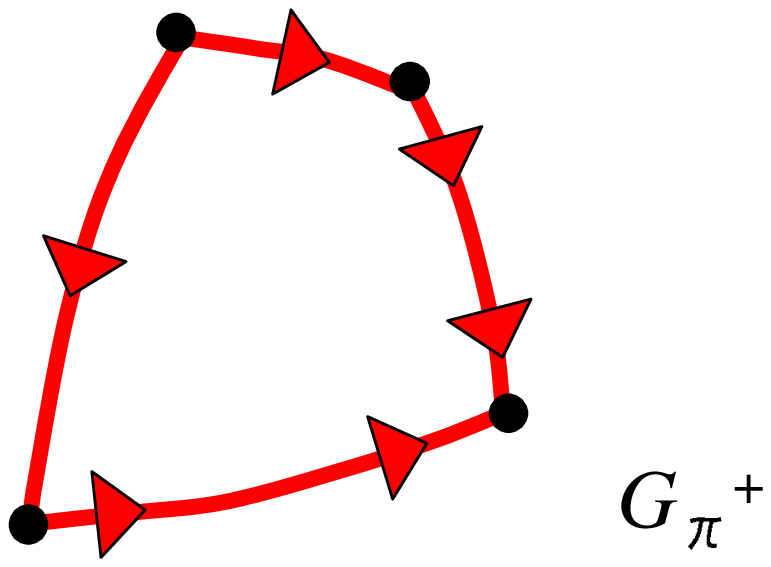}
\end{center}
\caption{An oriented matroid program $\pi = (M, g, f)$
and its OMP digraph $G_\pi^+$}
\label{fig:HKproperty}
\end{figure}
\\
If $M$ is representable,
its OMP digraph for any choice of $g$ and $f$ satisfies
the three properties of LP digraphs~\cite{om, EF82},
(1) acyclicity, (2) unique sink\&source property~\cite{SW01},
and (3) the Holt-Klee condition~\cite{HK},
but the situation is different in general.
Every OMP digraph satisfies the unique sink\&source property~\cite[page 426]{om},
but there exists an OMP digraph with a directed cycle,
see non-BOMs ~\cite{EF82} and non-Euclidean OMs \cite{om}.
On the other hand,
it is not known 
whether every OMP digraph $G_\pi^+$ satisfies
the Holt-Klee condition.
This motivates us to define the HK property for oriented matroids:
\begin{moriyamadf}
\label{df:HK-property_pi}
A proper oriented matroid program $\pi = (M, g, f)$ is called {\em HK}
if $G_\pi^+$ satisfies the Holt-Klee condition where
the dimension $d$ is $r(M)-1$,
and {\em non-HK} otherwise.
\end{moriyamadf}
\begin{moriyamadf}
\label{df:HK-property}
An oriented matroid $M = (E, {\mathcal L})$ is called {\em HK}
if the oriented matroid program $\pi= (M, g, f)$ is HK
for any two distinct elements $f, g \in E$ for which
$\pi$ is proper
and the same holds for any reorientations and any minors of $M$,
and {\em non-HK} otherwise.
\end{moriyamadf}
As we mentioned above,
if $M$ is representable,
every OMP digraph satisfies the Holt-Klee condition,
hence we obtain the following proposition:
\begin{prp}
Every representable oriented matroid has the HK property.
\end{prp}
For $r \leq 3$,
every unique-sink unique-source orientation on the graph of $P_\pi$
satisfies the Holt-Klee condition.
This observation gives the following proposition:
\begin{prp}
\label{prp:NOnon-HK}
Every oriented matroid of rank $r$ has the HK property
if $r \leq 3$.
\end{prp}

Now we consider the dual interpretation of the HK
condition which leads to the notion of HK* property.  
For this, we use the facet graph of a convex
polytope with orientation induced
by a line shelling ordering
of facets, see \cite{BM71} and \cite[Theorem 8.11]{Ziegler}.  
In the setting of oriented matroids,
the role of a straight line in an arrangement
of hyperplanes can be played by a coline.  What
we shall obtain is a \emph{coline shelling},
which is a special kind of what is known as tope graph
shelling or pseudoline shelling \cite{EF82}, see 
\cite[Section 4.3]{om} for a more algebraic treatment.

A \emph{coline fixation} is a pair $\omega = (M, T)$,
where $M$ is an oriented matroid $M=(E, \LL)$
and $T \subset E$ is a coline of $M$.
The associated \emph{supercell} $P_\omega$ of $(M,T)$ is
$P_\omega := \{ X \in \LL :
X_e = \{ 0, + \}~\textrm{for all}~e \in E \setminus T \}$.
A vector $Z \in \LL$ is an \emph{interior point} of $P_\omega$
if $Z^+ = E \setminus T$.
For each element $f \in E \setminus T$, the subset of
$P_\omega$ defined by
$P_\omega(f):=\{ X \in P_\omega :
X_f=0,
X^+ = E \setminus \{ T \cup \{ f \} \}$ is the \emph{face}
of $P_\omega$ induced by $f$.  The facets of $P_\omega$ are
the faces $P_\omega(f)$ that are maximal.  We say
$T$ is \emph{generic} in $M$ if
there exists $Y \in \LL$ for every $f \in E \setminus T$ such that
$Y^0 = T \cup \{ f \}$ and
$\underline{Y} = E \setminus (T \cup \{ f \})$.
A coline fixation $\omega = (M, T)$ is called \emph{proper}
if $T$ is generic, there exists an interior point 
$Z \in P_\omega$ such that $Z^0 = T$, and 
all faces $P_\omega(f)$ ($f\in E \setminus T$)
are facets.

When a coline fixation $\omega = (M,T)$ is proper,
there is a unique linear ordering of the elements
of $E \setminus T$
up to reversal:
$e_1, e_2, ..., e_s$ ($s={|E\setminus T|}$)
such that for each $k=1,\dots,s$,
the vector $V^k$ defined by
$(V^k)^- = \{ e_i : 1 \leq i \leq k-1  \}$
and $(V^k)^0 = T \cup \{ e_k \}$ is a cocircuit of $M$.
This ordering (unique up to reversal) 
is called the \emph{coline shelling} 
induced by $\omega$, denoted by $CS_\omega$.
By the duality of the ranking of vertices in a convex
polytope and the line shelling of the dual polytope,
we  define the facet graph with orientation
induced by a shelling.
Namely, we define 
the \emph{shelling digraph} $SG_\omega$ of $\omega$ as follows:
The set of vertices of $SG_\omega$ is $E \setminus T=\{e_1, \dots, e_s\}$,
and there is an edge $(e_i, e_j)$ directed from $e_i$ to $e_j$
if and only if $i<j$ and the associated two facets 
$P_\omega(e_i)$ and $P_\omega(e_j)$ are
\emph{adjacent}, i.e., their intersection is maximal over
all intersections of two distinct facets.
\begin{figure}[ht]
\begin{center}
\includegraphics[bb=70 400 420 723,scale=0.38]{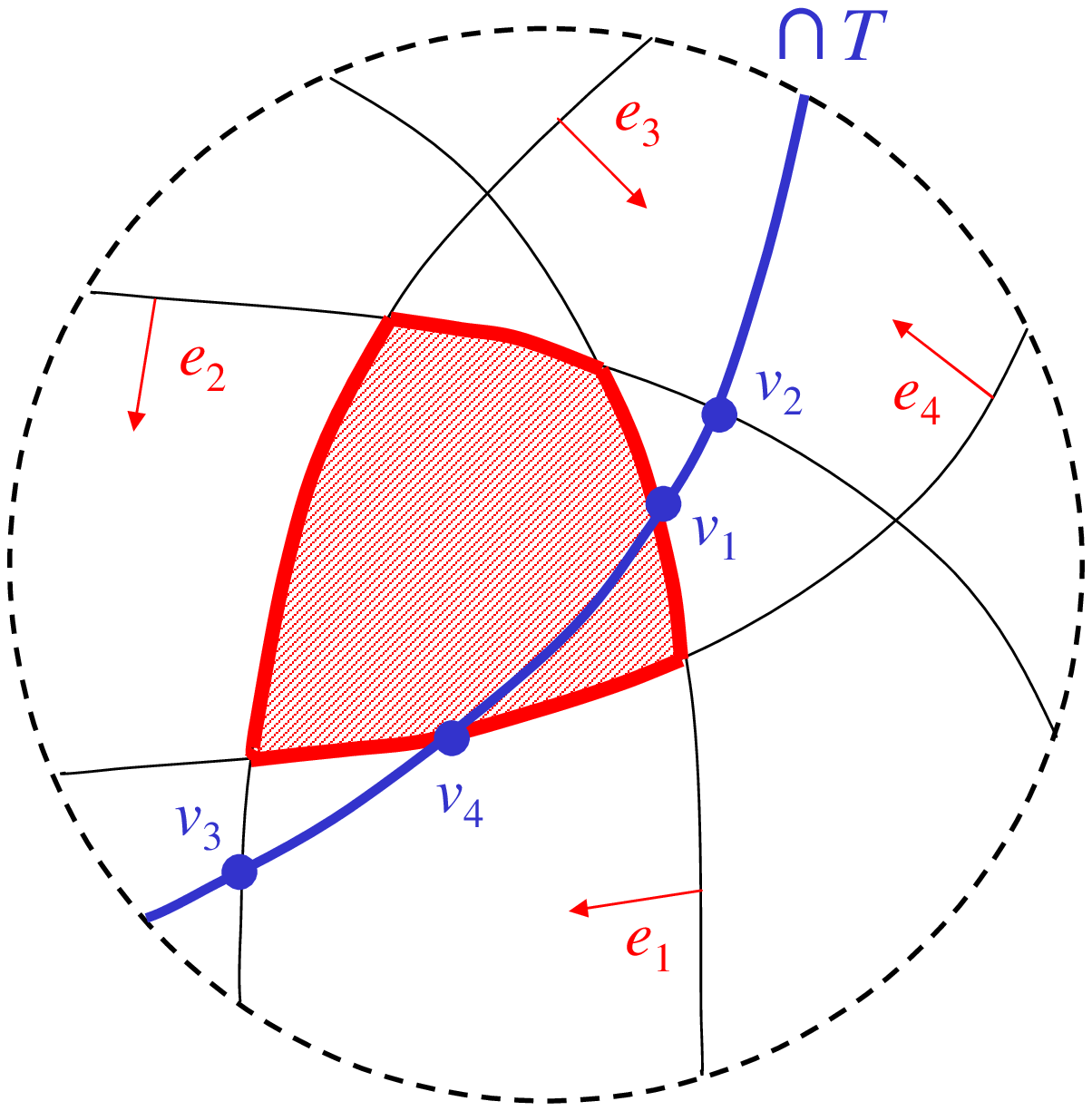}
\hspace{8mm}
\includegraphics[bb=120 400 362 723,scale=0.38]{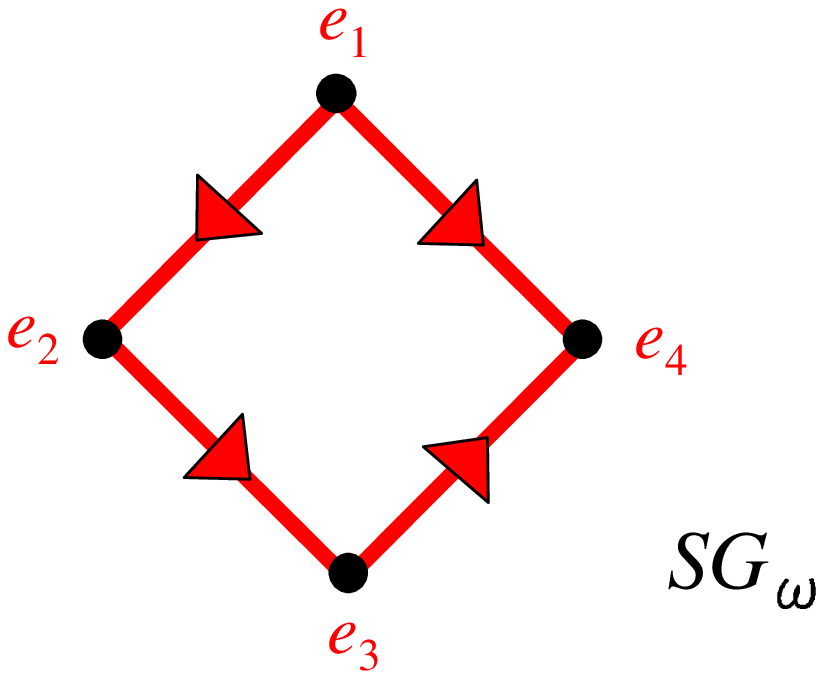}
\end{center}
\caption{A coline fixation $\omega = (M, T)$
and its shelling digraph $SG_\omega$}
\label{fig:dualHKproperty}
\end{figure}
\\
If $M$ is representable, the arrangement is
realizable as a hyperplane arrangement, and thus
every coline shelling $CS_\omega$ is realizable
as a line shelling.
By duality, the shelling digraph is an LP digraph~\cite{Ziegler}.
Hence every shelling digraph of a proper coline fixation
in a representable oriented matroid
satisfies (1) acyclicity, (2) unique sink\&source property~\cite{SW01}, and
(3) the Holt-Klee condition~\cite{HK}.
It is worthwhile to observe that even 
when $M$ is non-representable,
every shelling digraph satisfies
acyclicity and the unique sink\&source property \cite{M06}.
We will show that not every shelling digraph $SG_\omega$ satisfies
the Holt-Klee condition.
To make our claim clear, it is important to
define the HK* property of an oriented matroid $M$:
\begin{moriyamadf}
\label{df:HK*-property_omega}
A proper coline fixation $\omega = (M, T)$ is called {\em HK*}
if $SG_\omega$ satisfies the Holt-Klee condition,
and {\em non-HK*} otherwise.
\end{moriyamadf}
\begin{moriyamadf}
\label{df:HK*-property}
An oriented matroid $M = (E, {\mathcal L})$ is called {\em HK*}
if the coline fixation $\omega= (M, T)$ is HK*
for any coline $T \subset E$ such that 
$\omega$ is proper,
and the same holds for any reorientations of $M$ and any minors of $M$,
and {\em non-HK*} otherwise.
\end{moriyamadf}
As we mentioned above,
if $M$ is representable,
every shelling digraph satisfies the Holt-Klee condition,
hence we obtain the following proposition:
\begin{prp}
Every representable oriented matroid has the HK* property.
\end{prp}
Likewise Proposition \ref{prp:NOnon-HK},
we also obtain the following proposition:
\begin{prp}
\label{prp:NOnon-HK*}
Every oriented matroid of rank $r$ has the HK* property
if $r \leq 3$.
\end{prp}

By Definition \ref{df:HK-property} and Definition \ref{df:HK*-property},
every representable oriented matroid belongs to both
the class of HK oriented matroids and the class of HK* oriented matroids.
In other words, all non-HK and all non-HK* oriented matroids are non-representable.
However, a non-representable oriented matroid is not necessarily
non-HK or non-HK*.
To understand these two classes more clearly, 
we look at the oriented matroids of rank $4$
on an $8$-element ground set.
We denote by OM(4,8) the class of all (non-isomorphic)
oriented matroids of rank $4$ on an $8$-element ground set.
The class OM(4,8) is the smallest (with respect to the rank and at the same time
the size of a ground set) that contain a non-representable oriented matroid.
Finschi and Fukuda gave a complete enumeration of the oriented matroids
in OM(4, 8) including non-uniform oriented matroids,
and we utilize their list~\cite{FF01,FF03,Fin_web}.
Note that an oriented matroid or rank $r$ is called 
\emph{uniform}
if every subset of cardinality $r$ is a basis,
otherwise \emph{non-uniform}.

The class OM(4,8) contains $2,628$ uniform oriented matroids,
and $24$ out of them are non-representable~\cite{R88,BRG90b}.
By a computer program we found
$18$ non-HK* oriented matroids
while there exist no non-HK oriented matroids.
Furthermore, the class OM(4,8) contains $178,844$ non-uniform oriented matroids,
and we found that $1,364$ out of them have the non-HK* property
and none has the non-HK property.

Finally, we show how one can use \emph{sensitive LP digraphs}
to construct an infinite family of
non-HK* oriented matroids.
The word ``sensitive'' means that
some minor change makes the LP digraph violate the Holt-Klee condition.
It is interesting to note, however, that
we have not found any non-HK oriented matroid so far.
If every oriented matroid is HK,
then it implies that
the Holt-Klee theorem may be provable in a 
purely combinatorial manner using
the oriented matroid axioms only.
We leave this question as an open problem.

\section{Enumeration of non-HK and non-HK* oriented matroids}

In this section,
we give a classification of oriented matroids
on an $8$-element ground set $E$ of rank $4$ in terms
of HK and HK* properties.  We also 
compare these two properties
with the existing properties of representable 
oriented matroids.

For the classification,
we use the database of oriented matroids 
by Finschi and Fukuda~\cite{FF01,FF03}
that contains both 
all uniform and non-uniform oriented matroids.
Table \ref{FFdb} shows the number of non-isomorphic oriented matroids.
It enumerates the oriented
matroids up to isomorphism of the associated big face lattices
(see \cite[Section 4.1]{om}).  Thus, in particular,
two oriented matroids equivalent by a reorientation
or by a permutation of the ground set are isomorphic.
Since the HK and the HK* properties as well as
representability are closed under such
operations, the database is well suited for our purpose.

Denote by OM($r$, $n$) the class of
non-isomorphic oriented matroids of rank $r$
on an $n$-element ground set.
From Proposition \ref{prp:NOnon-HK} and Proposition \ref{prp:NOnon-HK*},
if $r\leq 3$, then OM($r$, $n$) contains no non-HK or non-HK* oriented matroid.
Thus, to seek a non-HK or a non-HK* oriented matroid, the rank $r$ has to be
at least four.
Since every rank-$4$ oriented matroid is representable if 
$n\leq 7$~\cite[Corollary 8.3.3]{om},
OM(4,8) is
the first candidate class that may contain a non-HK or non-HK* oriented matroid.
In $181,472$ oriented matroids of OM(4,8),
the number of uniform ones is $2,628$
and the number of non-uniform ones is $178,844$.
Bokowski and Richter-Gebert showed that
there exist $24$ non-representable uniform oriented matroids
among the $2,628$ uniform oriented matroids~\cite{BRG90b}.
Nakayama, Moriyama, Fukuda and Okamoto reconfirmed their
result using biquadratic final polynomials with the rational arithmetic~\cite{NMFO05}.
On the other hand,
the number of all non-representable non-uniform oriented matroids
is not known.

{\small
\label{FFdb}
\begin{table}
\begin{center}
\begin{tabular}{l| r r r r r r r r r r}
$|E|=$ & $1$ & $2$ & $3$ & $4$ & $5$ & $6$ & $7$ & $8$ & $9$ & $10$\\ \hline \hline
$r=1$ & $1$ & $$ & $$ & $$ & $$ & $$ & $$ & $$ & $$ & $$ \\
$r=2$ & $$ & $1$ & $1$ & $1$ & $1$ & $1$ & $1$ & $1$ & $1$ & $1$ \\
$r=3$ & $$ & $$ & $1$ & $2$ & $4$ & $17$ & $143$ & $4890$ & $461053$ & $95052532$ \\
$r=4$ & $$ & $$ & $$ & $1$ & $3$ & $12$ & $206$ & $\textbf{181472}$ & $$ & $$ \\
$r=5$ & $$ & $$ & $$ & $$ & $1$ & $4$ & $25$ & $6029$ & $$ & $$ \\
$r=6$ & $$ & $$ & $$ & $$ & $$ & $1$ & $5$ & $50$ & $508321$ & $$ \\
$r=7$ & $$ & $$ & $$ & $$ & $$ & $$ & $1$ & $6$ & $91$ & $$ \\
$r=8$ & $$ & $$ & $$ & $$ & $$ & $$ & $$ & $1$ & $7$ & $164$ \\
$r=9$ & $$ & $$ & $$ & $$ & $$ & $$ & $$ & $$ & $1$ & $8$ \\
$r=10$ & $$ & $$ & $$ & $$ & $$ & $$ & $$ & $$ & $$ & $1$ \\
\end{tabular}
\caption{
The number of non-isomorphic oriented matroids on a ground set $E$ of rank $r$~\cite{FF01,FF03}}
\end{center}
\end{table}
}

First, by our computer program, we enumerate
 the non-HK and the non-HK* oriented matroids in OM(4,8).
We found $18$ non-HK* uniform oriented matroids
out of $2,628$ uniform oriented matroids,
and $1,364$ non-HK* non-uniform oriented matroids
out of $178,844$ non-uniform oriented matroids.
On the other hand,
there exist no non-HK oriented matroids in OM(4,8).

Here, we present one non-HK* uniform oriented matroid IC(8,4,2).
Note that IC(n, r, c) refers to the $c$-th oriented matroid
of rank $r$ on the ground set $E$ such that $|E|=n$
in the catalog by Finschi and Fukuda~\cite{FF03,Fin_web}.
This OM is known as RS(8), constructed by Roudneff and Sturmfels \cite{RS88}.
Table \ref{chirotope_2} shows the chirotope representation
 of IC(8,4,2).  For example, the second column of
 the table indicates that the sign of the basis $\{1,2,3,5\}$ is $+$.  It is clear from
 the table that the OM is uniform, because there is no quadruple 
 taking the zero sign.

 \begin{table}[H]
  \begin{center}
   {\small
   \vspace{0mm}
   \verb||\\
   \verb|1111211121121231112112123112123123411121121231121231234112123123412345|\\
   \vspace{-1.5mm}
   \verb|2223322332334442233233444233444555522332334442334445555233444555566666|\\
   \vspace{-1.5mm}
   \verb|3344434445555553444555555666666666634445555556666666666777777777777777|\\
   \vspace{-1.5mm}
   \verb|4555566666666667777777777777777777788888888888888888888888888888888888|\\
   \vspace{-1.0mm}
   \verb|+++++++++++++++++++++++++++-------++------------+--++-----+--+---+--++|\\
   \vspace{-4.5mm}
   }
  \end{center}
    \caption{The chirotopes of IC(8,4,2)}
  \label{chirotope_2}
 \end{table}
\noindent
For a coline fixation $\omega = (IC(8,4,2), \{ 1,8 \})$,
the coline shelling $CS_\omega$ is the sequence $3, 2, 7, 6, 4, 5$ (up to reversal).
The coline shelling yields the following shelling digraph $SG_\omega$.

\begin{figure}[H]
\begin{center}
\includegraphics[bb=120 400 362 723,scale=0.38]{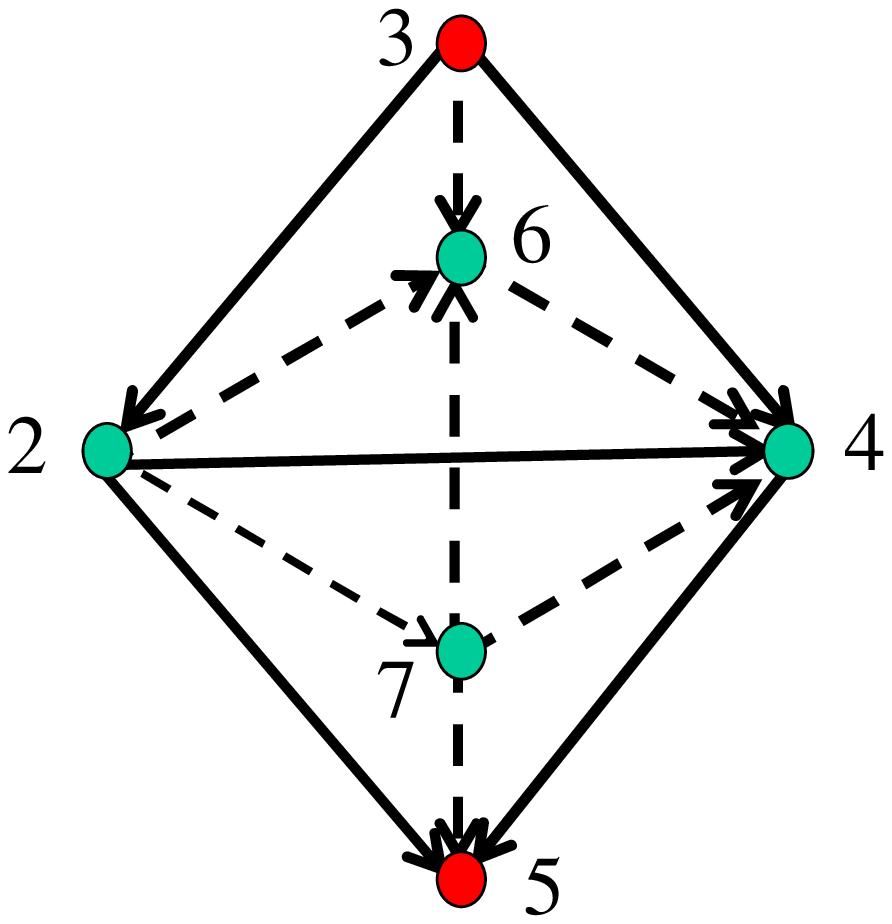}
\end{center}
\caption{The shelling digraph for a coline fixation $\omega = (IC(8,4,2), \{ 1,8 \})$}
\label{fig:RS8}
\end{figure}

\noindent
There are only two internally vertex-disjoint paths, since
all dipaths from $s$ to $t$ must go through at least one of the 
vertices $2$ or $4$.
Therefore, IC(8, 4, 2) is a non-HK* oriented matroid.

Secondly, we enumerate the oriented matroids with 
the two known properties of representable oriented
matroids,
the Euclidean \cite{EF82, EM82}
and Shannon properties~\cite{Sh76, Sh79},
in OM(4,8),
and compare them with the non-HK and non-HK* properties.
Let us recall the notion of Euclidean and Shannon 
matroids.

A oriented matroid program $\pi = (M, g, f)$ is called \emph{Euclidean}
if there exists no directed cycle in $G_\pi$,
and \emph{non-Euclidean} otherwise. An oriented matroid
$M = (E, {\mathcal L})$ is called \emph{Euclidean}
if $\pi= (M, g, f)$ is Euclidean for any two distinct elements $f \neq g \in E$.  By definition, the class is
closed under reorientations and taking minors.
Also, it is easy to see that 
every representable oriented matroid is Euclidean~\cite{EF82,EM82}.

The Shannon property is naturally defined
by a theorem by Shannon~\cite{Sh76,Sh79},
stating that every representable oriented matroid has
simplicial topes at least as many as twice the size of the ground set.

\begin{figure}[ht]
\begin{center}
\includegraphics[scale=0.40]{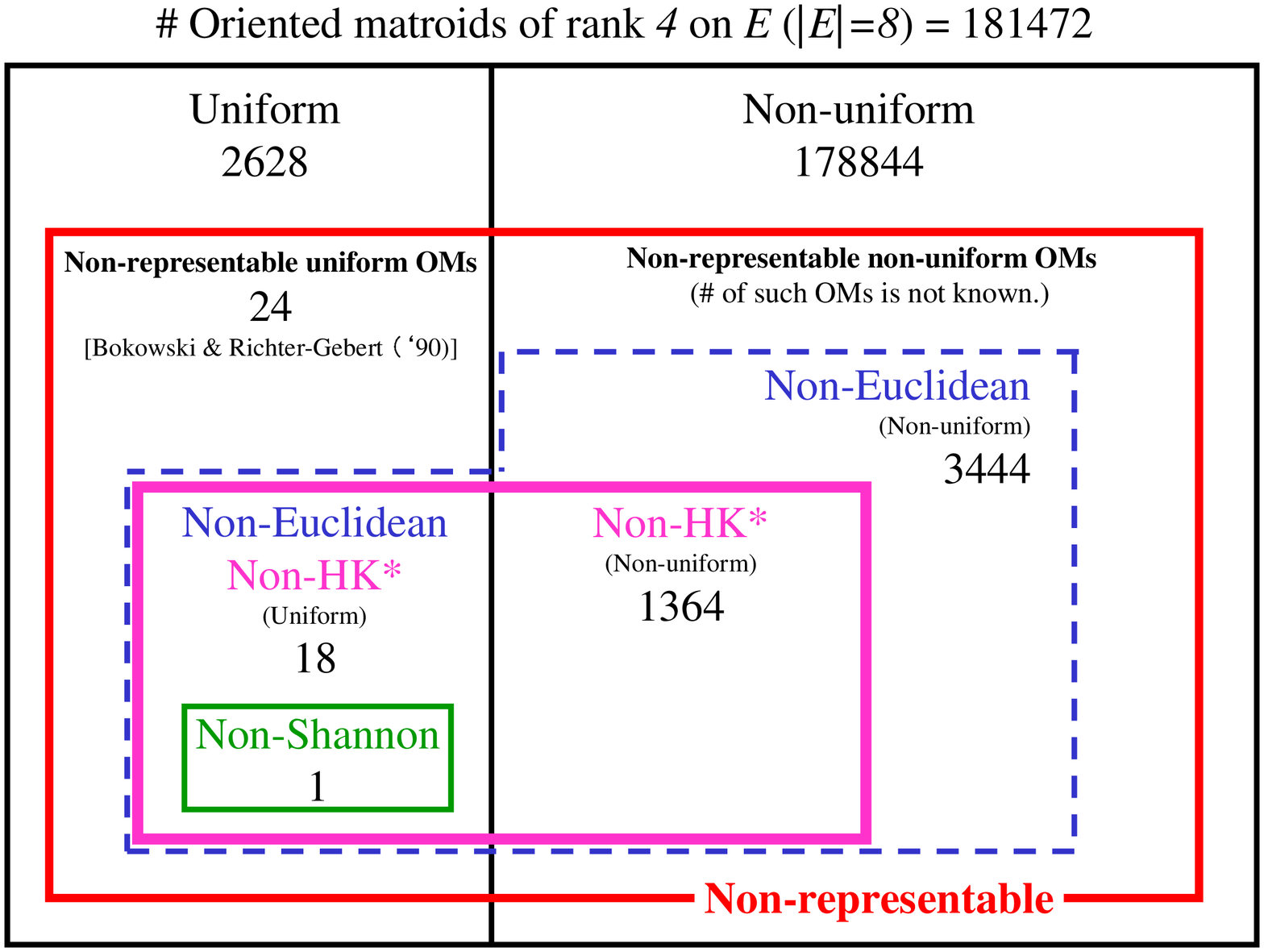}
\caption{A Classification of OM(4, 8)} \label{fig:experiment}
\end{center}
\end{figure}

We found $18$ non-Euclidean uniform oriented matroids,
and they are the same as the $18$ non-HK* uniform oriented matroids.
Furthermore, we found $3,444$ non-Euclidean non-uniform oriented matroids,
and they properly include all $1,344$ non-HK* non-uniform oriented matroids.
On the other hand, there exists
only one non-Shannon oriented matroid in OM(4,8),
which is also non-HK* and furthermore non-Euclidean.
This is known as RS(8) constructed by Roudneff and Sturmfels 
in \cite{RS88}.
Figure \ref{fig:experiment} summarizes the results.
For more precise information, we suggest the reader to
look at the web site 
\begin{quote}
\htmladdnormallink
{http://www-imai.is.s.u-tokyo.ac.jp/\~{}nak-den/OMcatalog/index.html}
  {http://www-imai.is.s.u-tokyo.ac.jp/\~{}nak-den/OMcatalog/index.html}
\end{quote}
which maintains the best of our knowledge.

\section{Construction of an infinite family of non-HK* oriented matroids}

In this section, we prove the following theorem.
\begin{moriyamathm}
\label{thm:mainInfiniteFamily}
For every $r \geq 4$ and every $n \geq 2r$,
there exists a non-HK* oriented matroid of rank $r$
on the ground set $E$ such that $|E| = n$.
\end{moriyamathm}
The essential idea of this proof is the notion
of sensitive LP digraphs,
which are special LP digraphs that
can lose the Holt-Klee property only by one flip.

Let $P$ be a $d$-polytope in $\R^d$,
$f(\vect{x}) = \vect{c}^T \vect{x}$\; a generic objective function,
$s$ the vertex of $P$ attaining the smallest value of $f$,
and $w$ the vertex attaining the second smallest value.
Notice that $( s, w )$ is an edge of $P$.
The quadruple $\gamma=(P, f, s, w)$ represents the 
LP digraph with two special vertices marked, which
will be called a \emph{marked} LP digraph.

\begin{moriyamadf}
\label{df:sensitiveLP}
A marked LP digraph $\gamma=(P, f, s, w)$ 
is called a \emph{sensitive} LP digraph
if by reversing (flipping) the orientation of 
the edge $(s, w)$, the resulting digraph violates
the Holt-Klee condition.
\end{moriyamadf}
We observe that both
 (1) acyclicity and (2) unique sink\&source property~\cite{SW01} 
remain satisfied
after reversing the orientation of the edge $(s, w)$.
A digraph satisfying the unique sink\&source property is called
a \emph{USO\/} digraph.
In the case of $d=2$,
any acyclic USO digraph satisfies the Holt-Klee condition.
Hence there exist no sensitive LP digraphs.
On the other hand, in $d \geq 3$,
there exists an acyclic USO digraph
not satisfying the Holt-Klee condition,
as in Figure \ref{fig:HK-non-HK}.
Based on the enumeration of combinatorial types of polytopes
with respect to the dimension and the number of vertices
by Finschi and Fukuda~\cite{FF01,FF03},
we have checked whether a $3$-polytope 
with few vertices admits a sensitive LP digraph.
When the number of vertices is less than six,
all acyclic USO digraphs satisfy the Holt-Klee condition.
\begin{figure}[ht]
\begin{center}
\includegraphics[scale=0.27]{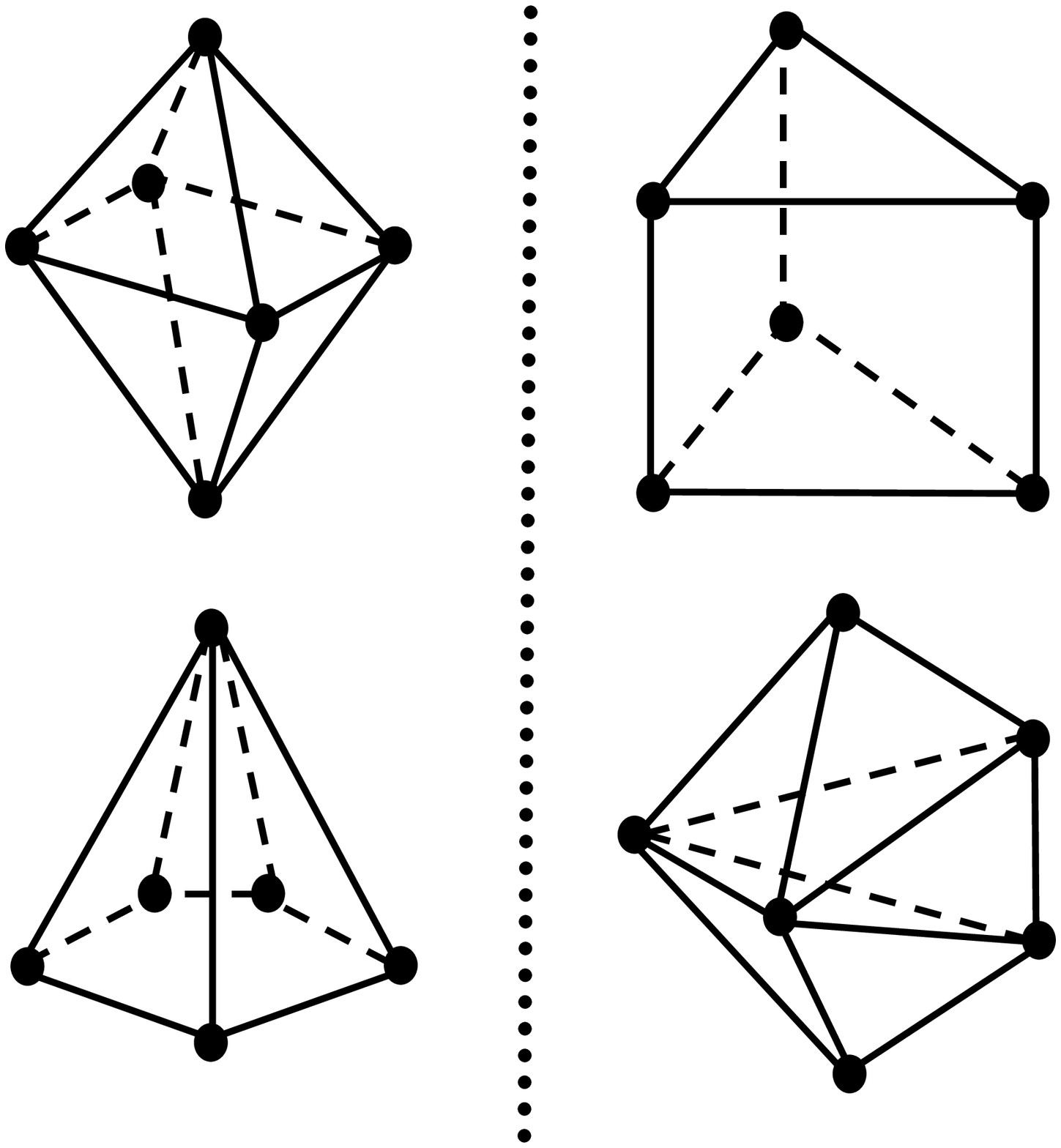}
\includegraphics[scale=0.27]{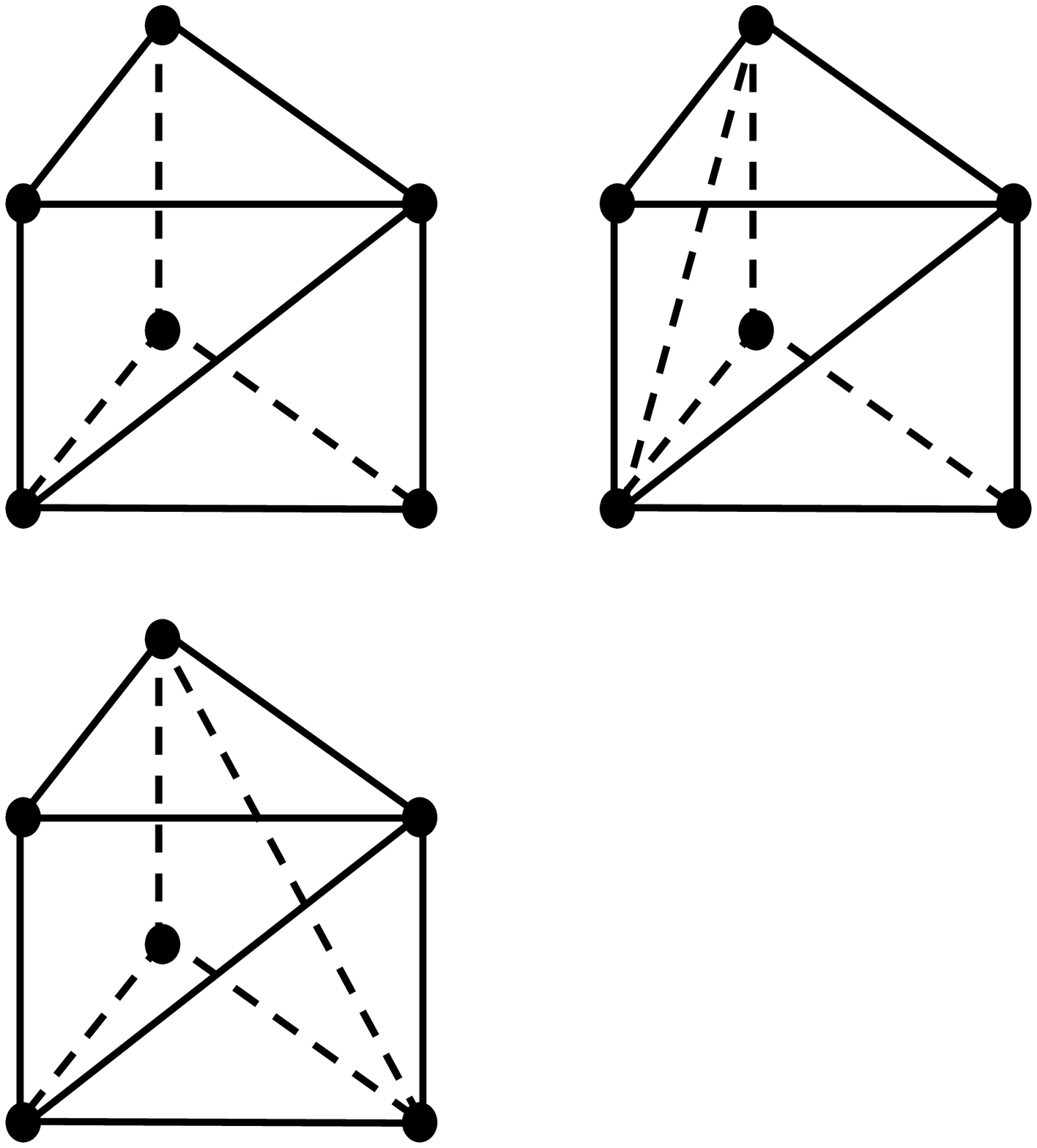}
\end{center}
\caption{Seven types of $3$-polytopes with six vertices.}
\label{fig:seven}
\end{figure}

\begin{figure}[ht]
\begin{center}
\includegraphics[scale=0.45]{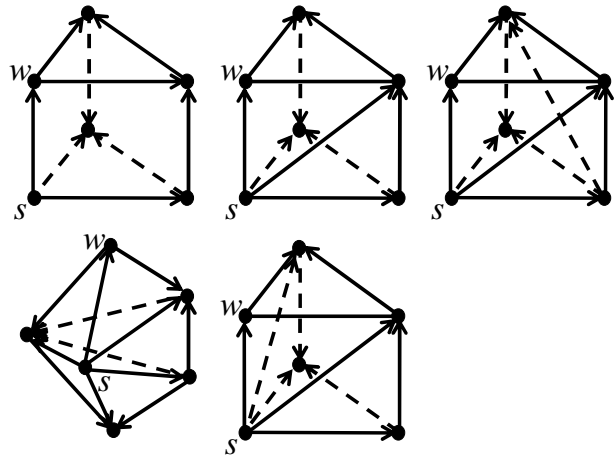}
\end{center}
\caption{Smallest sensitive LP orientations
w.r.t the dimension and the number of vertices.
(The leftmost two orientations appear in \cite{TI99} and \cite{Si90}.)}
\label{fig:seven2}
\end{figure}

However, things are different if the number of vertices is equal to six.
Among the seven types of $3$-polytopes with six vertices in Figure \ref{fig:seven},
all (five polytopes) except for the leftmost two polytopes admit
sensitive LP orientations, as shown 
in Figure \ref{fig:seven2}.
Thus, we have the following.
\begin{prp}
\label{prp:smallestSensitiveLP}
The sensitive LP digraphs $\gamma=(P, f, s, w)$ in Figure \ref{fig:seven2}
are smallest with respect to the dimension of $P$ and the number of vertices of $P$.
\end{prp}

By using sensitive LP digraphs,
we construct an infinite family of non-HK* oriented matroids
through the following three steps.
First, we prove Theorem \ref{thm:main_flipping},
which states that if a $d$-polytope with $n$ vertices
admits a sensitive LP digraph,
then there exists a non-HK* oriented matroid of rank $r = d{+}1$
on a $(d{+}n{-}1)$-element ground set.
Second, we prove in Lemma \ref{lem:truncation}
that if a $3$-polytope $P$ with $n$ vertices
admits a sensitive LP digraph,
then a $3$-polytope with $n{+}1$ vertices induced by $P$
by a \emph{truncation},
also admits a sensitive LP digraph.
Therefore, combining them with Proposition \ref{prp:smallestSensitiveLP},
we obtain Proposition \ref{prp:secondstep},
which states that there exists a $3$-polytope with $n$ vertices
admitting a sensitive LP orientation for every $n \geq 6$.
Third, we prove in Proposition \ref{prp:pyramid}
that if a $d$-polytope $P$ with $n$ vertices
admits a sensitive LP digraph,
a $(d{+}1)$-polytope with $n{+}1$ vertices obtained from $P$
by a \emph{pyramid\/} construction
also admits a sensitive LP digraph.
Therefore, combining them all,
we obtain Proposition \ref{prp:thirdstep}
stating that there exists a $d$-polytope with $n$ vertices
admitting a sensitive LP digraph for every $d \geq 3$ and $n \geq d{+}3$.
Finally, we derive Theorem \ref{thm:mainInfiniteFamily}
from Theorem \ref{thm:main_flipping} and Proposition \ref{prp:thirdstep}.

\subsection{Construction of non-HK* OMs from representable OMs by a flipping}
\label{firststep}

In this section, we prove the following theorem:
\begin{moriyamathm}
\label{thm:main_flipping}
Suppose there is a $d$-polytope with $n$ vertices
whose graph admits a sensitive LP orientation. Then
there exists a non-HK* oriented matroid of rank $r = d{+}1$
on the ground set of size $|E|= n{+}d{-}1$.
\begin{Proof}
Let $P$ be a $d$-polytope in $\R^d$,
$F_i$'s ($1 \leq i \leq n$) be the facets of $P$,
and $H_i$ be
the facet-supporting hyperplane of $F_i$.
Here we take a line $L$ in general position 
through the interior of $P$.  

First, we show that
for any two adjacent facets $F_a$ and $F_b$,
there exists a a $d$-simplex in $\R^d$ such that
$V_a = L \cap H_a$ and $V_b = L \cap H_b$ are its vertices,
and the other $d{-}1$ vertices are
on the relative interior of $F_a \cap F_b$, see Figure \ref{fig:flip}.
In fact,
we can take an arbitrary $(d{-}2)$-simplex $\Delta^{d{-}2}$
in the relative interior of $F_a \cap F_b$.
The vertex $V_a = L \cap H_a$ is not on the $(d{-}2)$-dimensional 
flat $H_a \cap H_b$, and
hence the convex hull of $V_a$ and $\Delta^{d{-}2}$ is 
a $(d-1)$-simplex contained in $H_a$.
Similarly, the vertex $V_b = L \cap H_b$ is not on $H_a$,
and thus the convex hull of $V_b$, $V_a$ and $\Delta^{d{-}2}$
is a $d$-simplex, which we denote by $\Delta^d$.
Because the $d{+}1$ vertices of $\Delta^d$ are affinely independent,
for every facet $G_j$ for $1 \leq j \leq d-1$ of $\Delta^{d{-}2}$,
the $d{-}1$ vertices of $G_j$, $V_a$ and $V_b$ are also affinely independent in $\R^d$.
Here, we denote by $T_j$ the hyperplane
determined by the $d{-}1$ vertices of $G_j$, $V_a$ and $V_b$.
Therefore, $H_a$, $H_b$ and $\{ T_j : 1 \leq j \leq d-1 \}$
are the supporting hyperplanes of the $d$-simplex $\Delta^{d}$,
and $L$ is the intersection of $\{ T_j : 1 \leq j \leq d-1 \}$.

Now, suppose there is a $d$-polytope $Q$ with $n$ vertices
admitting a sensitive LP orientation.
By duality, this implies that
we may suppose that the $d$-polytope $P$ 
above is a dual to $Q$ and that the line $L$ induces a shelling digraph
isomorphic to the sensitive LP digraph.
Particularly, we may suppose that $F_a$ and $F_b$ are the first
(or the last) two facets of 
the line shelling of $P$ induced by $L$,
see Figure~\ref{fig:flip}.  This implies that
each of the hyperplanes $\{ H_i : 1 \leq i \leq n \} \setminus \{ H_a, H_b \}$
does not intersect with the $d$-simplex $\Delta^d$. 

\begin{figure}[ht]
\begin{center}
\includegraphics[scale=0.40]{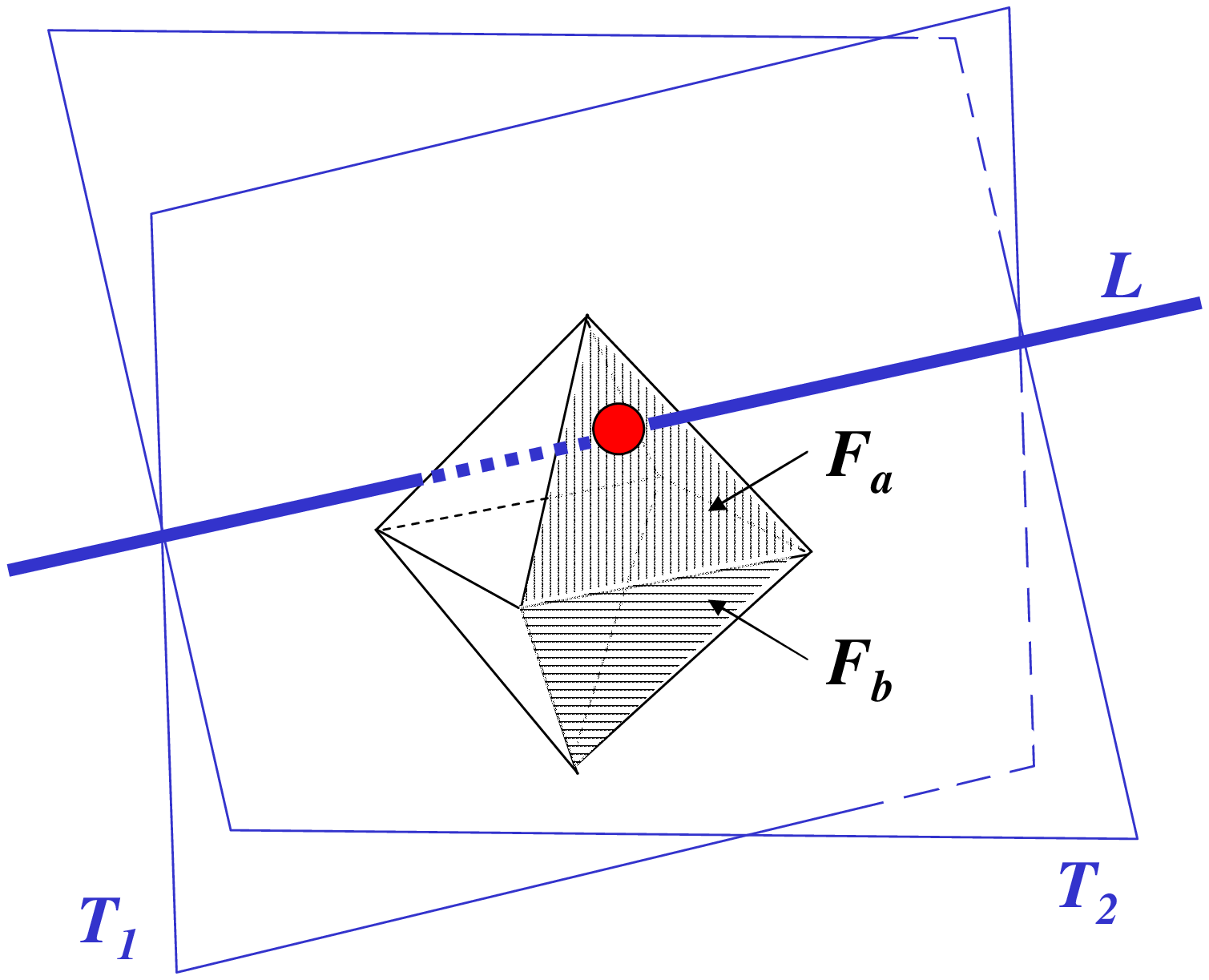}
\includegraphics[scale=0.40]{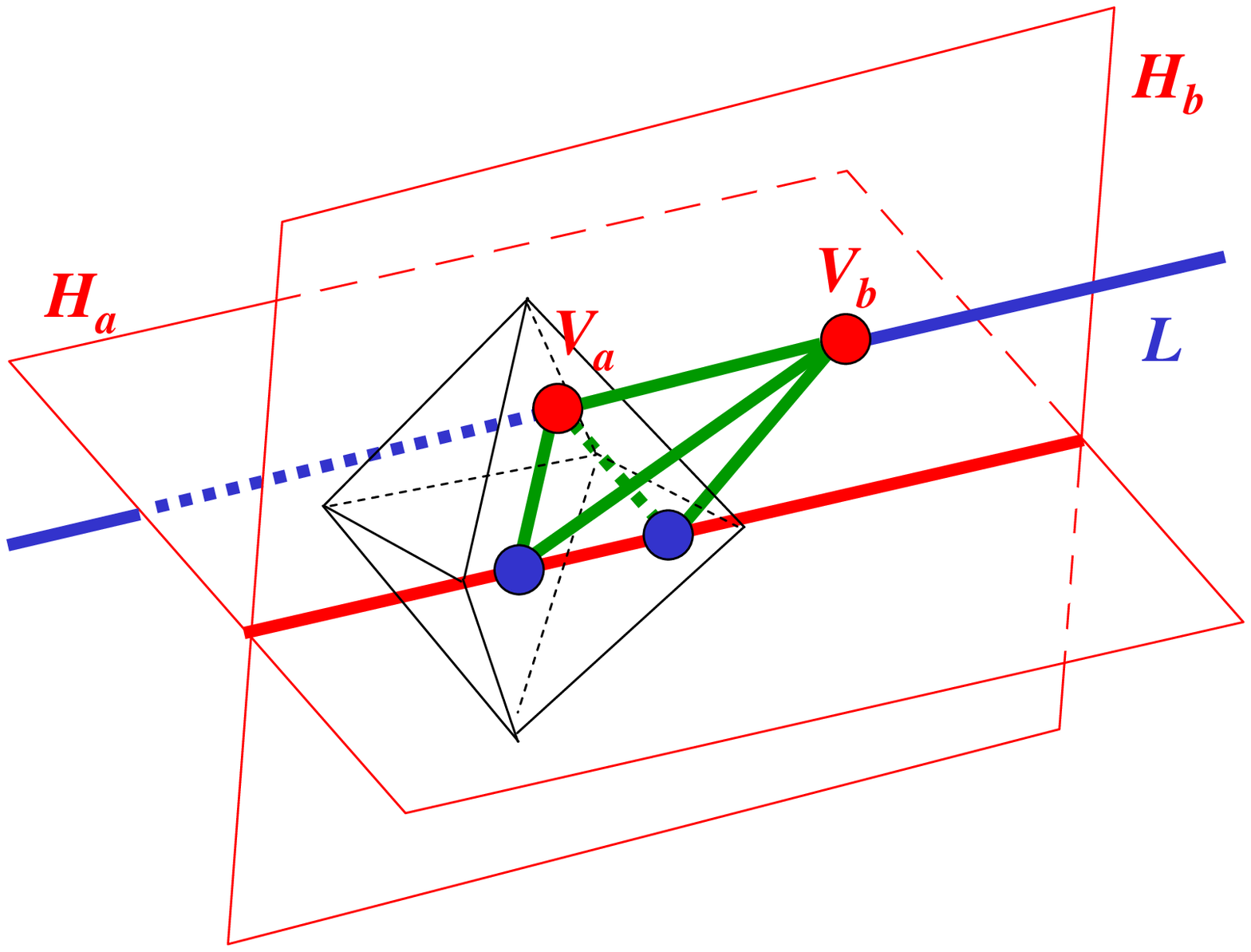}
\caption{The $d$-simplex for the first two facets of a line shelling} \label{fig:flip}
\end{center}
\label{fig:flipping}
\end{figure}

We claim that one can apply a flipping operation~\cite{FT88a}
to the representable oriented matroid $M = (E, {\mathcal L})$ of
the hyperplane arrangement ${\mathcal H}$ of
$\{ H_i : 1 \leq i \leq n \} \cup \{ T_j : 1 \leq j \leq d-1 \}$.
Let us present the precise definition of $M$.
For every $d$-dimensional hyperplane $h \in {\mathcal H}$,
we take a hyperplane $h'$ in $\R^{d+1}$
containing the origin and a lifted copy $H\times \{1\}$ of $H$.
$M$ is the representable oriented matroid of rank $d{+}1$
represented by the linear sphere arrangement
on the $d$-dimensional unit sphere $S^{d+1}$ in $\R^{d+1}$
of $\{ e_{i} = S^{d{+}1} \cap H_i' : 1 \leq i \leq n \}
\cup \{ e_{n{+}j} = S^{d{+}1} \cap T_j' : 1 \leq j \leq d-1 \}$,
i.e. $E = \{ e_k : 1 \leq k \leq n{+}d{-}1 \}$.
The subset $\{ e_{k} : n{+}1 \leq j \leq n{+}d{-}1 \} \subset E$ is a coline of $M$
by construction.
Let $\omega$ be the coline fixation $(M, T)$, where $T=\{ e_{n{+}j} : 1 \leq j \leq d-1 \}$. Then,
its coline shelling $CS_\omega$ coincides with the line shelling of $P$
given by $L$. In particular, the first two elements of $CS_\omega$ are $e_a$ and $e_b$.

Exploiting the structure of the hyperplane arrangement ${\mathcal H}$,
we apply a flipping operation~\cite{FT88a} to $M$.  Namely, we may
flip any element $e_k \in E' = \{e_a, e_b\} \cup \{ e_{n{+}j} : 1 \leq j \leq d{-}1 \}$
so that the associated simplex tope is flipped over.
The resulting oriented matroid does not depend on the choice
of $e_k$ and is denoted by $M'$.

Now, we observe that $T$ remains a coline in $M'$
and the coline shelling $CS_{\omega'}$ induced by
the fixation $\omega'=(M',T)$ differs from $CS_\omega$
only for the ordering of $e_a, e_b$.
This means that the orientation of the edge $(e_a, e_b)$
in the shelling digraph $SG_\omega$ is reversed in
$SG_{\omega'}$, and thus
$SG_{\omega'}$ does not satisfy the Holt-Klee condition. This means that
the oriented matroid $M'$ is non-HK*.  It has rank $r = d{+}1$
and $|E|= n{+}d{-}1$.  This completes the proof.
\end{Proof}
\end{moriyamathm}

Using both Proposition \ref{prp:smallestSensitiveLP}
and Theorem \ref{thm:main_flipping},
we have a theoretical proof for the fact we knew from
our computational classification.
\begin{cor}
There exists a non-HK* oriented matroid of rank $4$
on a $8$-element ground set.
\end{cor}

\subsection{Construction of non-HK* OMs of rank $4$ by a truncation}
\label{secondstep}

First, we define a truncated polytope, which is a key idea of this section.
\begin{moriyamadf}
\label{df:truncatedPolytope}
Let $P$ be a $3$-polytope in $\R^3$
containing a \emph{simple\/} vertex $v$ (i.e. a vertex $v$ with exactly $3$ neighbors),
$\{ v_i : i=1, 2, 3 \}$ the vertices adjacent to $v$
and $\{ u_j : j=1,2 \}$ points in the relative interior of an edge $(v, v_j)$.
By a \emph{truncated polytope\/} $tr(P)$, we mean a $3$-polytope 
$P \cap (H \cup H^+)$
where $H$ is the hyperplane determined by $u_1$, $u_2$ and $v_3$,
and $H^+$ is the open halfspace of $H$ containing all vertices except $v$,
see Figure~\ref{fig:truncation}.
\begin{figure}[ht]
\begin{center}
\includegraphics[scale=0.40]{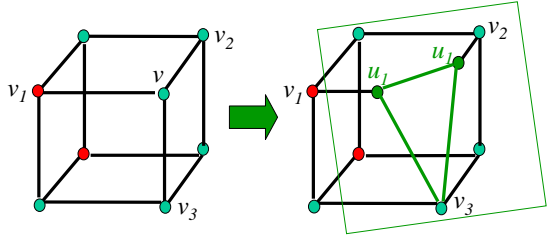}
\end{center}
\caption{A truncated polytope}
\label{fig:truncation}
\end{figure}
\end{moriyamadf}

\begin{moriyamalem}
\label{lem:truncation}
Let $P$ be a $3$-polytope in $\R^3$
containing a simple vertex $v$.
If $P$ admits a sensitive LP digraph,
a truncated polytope $tr(P)$ also admits a sensitive LP digraph.
\begin{Proof}
Let $(P, f, s, w)$ be a sensitive LP digraph.
A truncation operation with respect to a simple vertex $v$
generates the following five new edges:
$(v_1, u_1)$, $(v_2, u_2)$, $(v_3, u_1)$, $(v_3, u_2)$, and $(u_1, u_2)$.
We show how to orient the five edges so that $tr(P)$ also admits a sensitive LP digraph.
\begin{figure}[ht]
\begin{center}
\includegraphics[scale=0.35]{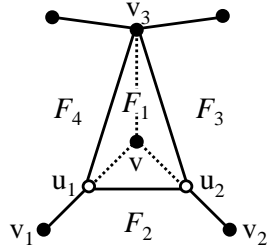}
\end{center}
\caption{The five new edges of $tr(P)$}
\label{fig:truncation_edges}
\end{figure}

By the symmetry of $v_1$ and $v_2$,
there are exactly 
six types of orientations of $(v, v_1)$, $(v, v_2)$, and $(v, v_3)$
with respect to outdegree and indegree of $v$.

In the case of (a), (b) and (c), since
no vertices of $v$, $v_1$, $v_2$ and $v_3$ are the global source $s$,
we only have to orient the five edges so that
(1) acyclicity, (2) the unique sink\&source property and (3) the Holt-Klee condition
are satisfied,
i.e., every facet of $F_1$, $F_2$, $F_3$ and $F_4$ in Figure \ref{fig:truncation_edges}
has a unique source and a unique sink,
and the number of disjoint paths between $v_1$, $v_2$ and $v_3$ remains unchanged.
Then, the resulting orientation of
$(v_1, u_1)$, $(v_2, u_2)$, $(v_3, u_1)$, $(v_3, u_2)$, and $(u_1, u_2)$
is classified into (a), (b) and (c)
in Figure \ref{fig:truncation_oriented_8cases}.
On the other hand, in the case of (d) and (e),
it is possible that
$v_2$ in (d) and $v_3$ in (e) are the global source $s$,
and $v$ in (d) and $v_3$ in (e) are $w$.
In these cases, if the five edges are oriented
as (d) and (e) in Figure \ref{fig:truncation_oriented_8cases},
i.e., $v_2$ in (d) and $v_3$ in (e) are also $s$,
and $u_2$ in (d) and $u_1$ in (e) are $w$,
the three properties are satisfied and
the sensitivity of an LP digraph remains unchanged.
Otherwise, as well as (a), (b) and (c),
we have only to orient the five edges such that
the three properties are satisfied.
Finally, in the case of (f),
$v$ is $s$ and one of $v_1$, $v_2$ and $v_3$ is $w$.
If the five edges are oriented
as (f) in Figure \ref{fig:truncation_oriented_8cases},
the three properties are also satisfied and
the sensitivity of an LP digraph also remains unchanged.
From the above, a truncated polytope $tr(P)$ also admits a sensitive LP digraph.
\begin{figure}[ht]
\begin{center}
\includegraphics[scale=0.48]{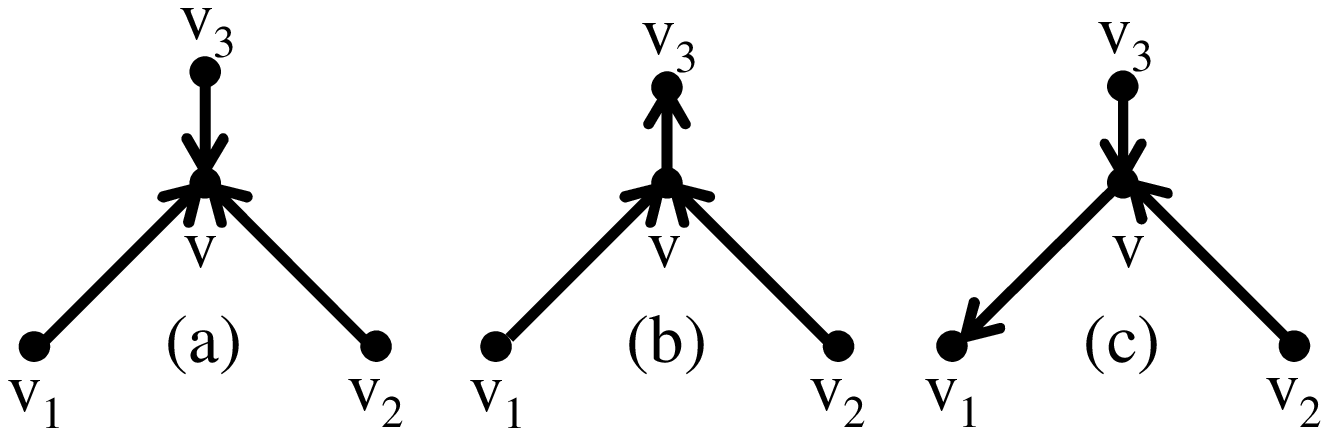}
\includegraphics[scale=0.448]{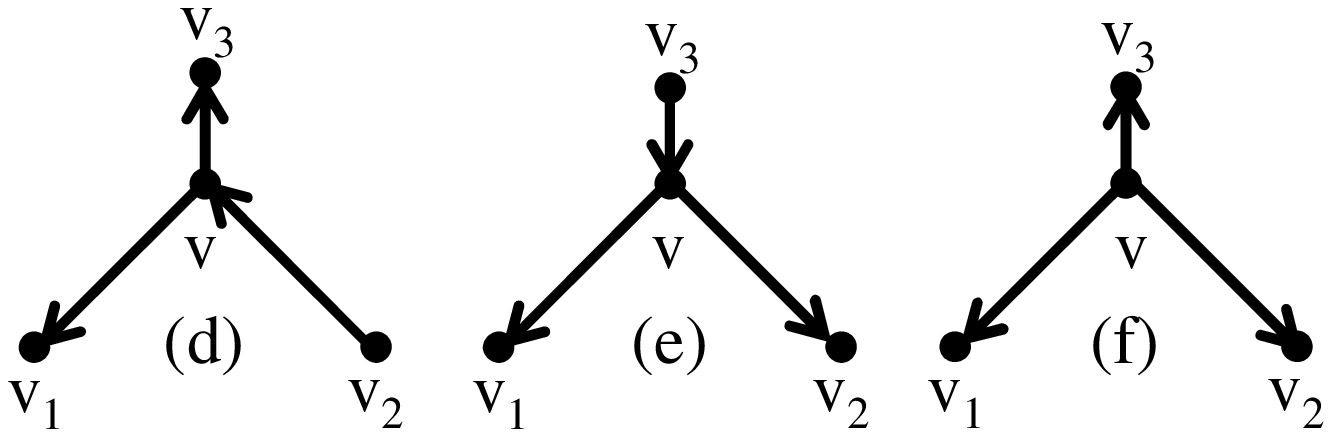}
\end{center}
\caption{The orientation of $(v, v_1)$, $(v, v_2)$, and $(v, v_3)$}
\label{fig:truncation_6cases}
\end{figure}
\begin{figure}[ht]
\begin{center}
\includegraphics[scale=0.55]{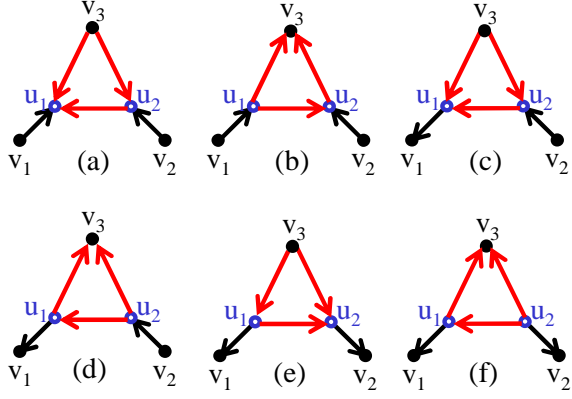}
\end{center}
\caption{The orientation of
$(v_1, u_1)$, $(v_2, u_2)$, $(v_3, u_1)$, $(v_3, u_2)$, and $(u_1, u_2)$}
\label{fig:truncation_oriented_8cases}
\end{figure}
\end{Proof}
\end{moriyamalem}
From Definition \ref{df:truncatedPolytope},
a truncation operation generates two simple vertices $u_1$ and $u_2$
while one simple vertex $v$ is removed.
Hence we have the following.
\begin{rem}
\label{rem:truncation}
There exists at least one simple vertex in a truncated polytope.
\end{rem}
Furthermore, because
all five polytopes in Proposition \ref{prp:smallestSensitiveLP}
contain a simple vertex,
we may apply a truncation operation to a $3$-polytope successively.
Thus, we obtain the main proposition of this section.
\begin{prp}
\label{prp:secondstep}
There exists a $3$-polytope with $n$ vertices which
admits a sensitive LP digraph for every $n \geq 6$.
\end{prp}

\subsection{Construction of non-HK* OMs of higher ranks}
\label{thirdstep}

In this section, we present a construction of a sensitive LP digraph
starting from a sensitive LP digraph in one lower dimension.

\begin{figure}[ht]
\begin{center}
\includegraphics[scale=0.35]{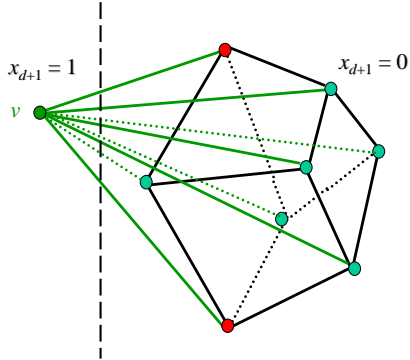}
\end{center}
\caption{An intuitive image of the pyramid polytope in $\R^4$}
\label{fig:pyramidPolytope}
\end{figure}

Given a $d$-polytope $P$ in $\R^d$,
its \emph{pyramid polytope\/} $py(P,v)$
is a $(d{+}1)$-polytope in $\R^{d+1}$
which is the convex hull of $P\times\{0\}$ and
a point $v\in \R^{d+1}$ not on the $d$-dimensional subspace containing $P$.
A canonical choice is to set $v_{d+1}=1$,
see Figure \ref{fig:pyramidPolytope}.

\begin{prp}
\label{prp:pyramid}
Let $P$ be a $d$-polytope in $\R^d$.
If $P$ admits a sensitive LP digraph,
a pyramid polytope $py(P, v)$,
also admits a sensitive LP digraph.
\begin{Proof}
Let $P$ be a $d$-polytope in $\R^d$ which admits a sensitive LP digraph.
This means $d\ge 3$ and thus $P$ has at least four vertices.
Let $\gamma=(P, f, s, w)$ be a sensitive LP digraph, where $f$ is a generic
objective function.
and let $z$ be the vertex of $P$ attaining the third smallest objective value.
Thus, $f(s) < f(w) < f(z)$.

Let $v\in \R^{d+1}$ be any point with $v_{d+1}=1$, and consider
the pyramid $Q = py(P,v)$.  We shall construct an objective function $g$
for $Q$ which induces a sensitive LP orientation.
Let $g$ be a natural extension of $f$: $g(\vect{y}) = f(\vect{x}) + c_{d{+}1}x_{d{+}1}$.  Since $v$ is the only vertex of $Q$ with nonzero last component,
one can set $c_{d{+}1}$ in such a way that
$g(S) < g(\vect{W}) < g(\vect{V}) < g(\vect{Z})$, where uppercase
letters denote the same (lowercase) vectors lifted to $\R^{d+1}$:
$\vect{S}^T = (\vect{s}^T, 0)$,
$\vect{W}^T = (\vect{w}^T, 0)$,
$\vect{Z}^T = (\vect{z}^T, 0)$,
$\vect{V}=v$, see Figure \ref{fig:pyramid_function}.

\begin{figure}[ht]
\begin{center}
\includegraphics[scale=0.35]{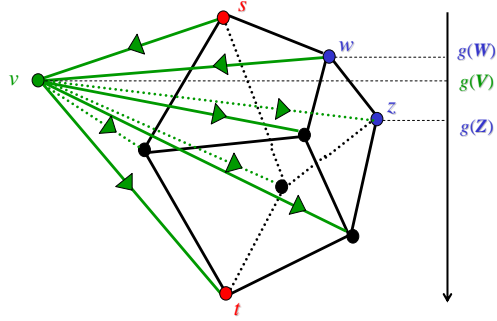}
\end{center}
\caption{A pyramid polytope and a generic function $g$}
\label{fig:pyramid_function}
\end{figure}

We claim that the LP digraph $Q$ induced by $g$ is sensitive,
or more precisely, the marked LP digraph $\gamma'= (Q, g, S, W)$
is a sensitive LP digraph.  By the construction, $S$,
$W$ and $T$ are the vertices attaining the smallest,
the second smallest, and the largest $g$ value, respectively.
The only difference between $\gamma$ and $\gamma'$ are
the extra edges in $\gamma'$ incident with $V$.  This means
that the maximum number of dipaths from $W$ to $T$ in $\gamma'$
with the edge $(S, W)$ reversed is at most one more than
the maximum number of dipaths from $w$ to $t$ in $\gamma$
with the edge $(s, w)$ reversed.  Since $\gamma$ is
sensitive, $\gamma'$ is sensitive as well.
\end{Proof}
\end{prp}
Combining with Proposition \ref{prp:secondstep},
we can construct a $d$-polytope which admits a sensitive LP digraph
successively.
Thus, we have the following.
\begin{prp}
\label{prp:thirdstep}
For each $d \geq 3$ and $n \geq d{+}3$,
there exists a $d$-polytope with $n$ vertices
which admits a sensitive LP digraph.
\end{prp}
Finally, from Theorem \ref{thm:main_flipping} and Proposition \ref{prp:thirdstep}
the main theorem of this section, Theorem \ref{thm:mainInfiniteFamily}, follows.

\section{Concluding remarks}

In this paper,
we introduced two new classes HK and HK* of oriented matroids
based on the Holt-Klee condition and its dual interpretation
in terms of line shellings, respectively.
In particular, the non-HK and non-HK* properties
are certificates for non-representability.

We have shown that these two classes are distinct.
While we gave a construction of an infinite family of non-HK* OMs
using the notion of sensitive LP digraphs,
it is not clear whether there exists any non-HK OM.
This leads to a fundamental question as to whether 
the Holt-Klee theorem can be proven combinatorially by using the OM axioms only
while the original proof in \cite{HK} relies heavily on geometric operations
such as affine transformations and orthogonal projections.

To get a better understanding,
we presented a classification of the oriented matroids of
rank $4$ on $8$-element ground set 
with respect to the HK, HK*, Euclidean and Shannon properties.
Our classification shows that there exists no non-HK OM in this class.
This suggests us to try to prove the statement
that every OM is HK.  A successful trial would 
yield a purely combinatorial proof of the Holt-Klee theorem.
\section*{Acknowledgments}
We wish to thank an anonymous referee for many constructive
suggestions.


\bibliographystyle{plain}
\bibliography{article}

\end{document}